\documentclass[11pt]{article}

\usepackage[font=small,labelfont=bf]{caption}
\usepackage{listings,spalign}
\usepackage{subcaption}
\usepackage{amsmath,amssymb}
\usepackage[left=1in, right=1in, top=1in, bottom=1in]{geometry}
\usepackage{verbatim}
\usepackage{titling}
\usepackage{parskip}% http://ctan.org/pkg/parskip
\usepackage{todonotes} % IZKISS

\usepackage{graphicx} % shane, figures
\usepackage[english]{babel}
\usepackage[utf8]{inputenc}
\usepackage{slashed}
\usepackage{amsfonts}% Caratteri Matematici
\usepackage{amsthm}% Gestione Teorem
\usepackage{amsmath,amssymb}
\usepackage{empheq}
\usepackage{dsfont}%Serve a fare \mathds{1} e ottenere l'identità
\usepackage{xcolor}
\usepackage{hyperref}
\usepackage{xcolor}

\usepackage{soul}

\newtheorem{prop}{Proposition}

\title{On parameter identifiability in network-based epidemic models}
\author{Istv\'an Z. Kiss$^1$ \& P\'eter L. Simon$^{2,3}$}
\date{$^1$ Department of Mathematics, University of Sussex, Falmer, Brighton BN1 9QH, UK\\
      $^2$ Institute of Mathematics, E\"otv\"os Lor\'and University Budapest, Hungary \\
$^3$ Numerical Analysis and Large Networks Research Group,\\ Hungarian Academy of Sciences, Hungary\\
      }

\begin{document}
\maketitle
%%%%%%%%%%%%%%%%%%%%%%%%%%%%%%%%%%%%%
\begin{abstract}
%%%%%%%%%%%%%%%%%%%%%%%%%%%%%%%%%%%%%
Many models in mathematical epidemiology are developed with the aim to provide a framework for parameter estimation and then prediction. It is well-known that parameters are not always uniquely identifiable. In this paper we consider network-based mean-field models and explore the problem of parameter identifiability when observations about an epidemic are available. Making use of the analytical tractability of most network-based mean-field models, e.g., explicit analytical expressions for leading eigenvalue and final epidemic size, we set up the parameter identifiability problem as finding the solution or solutions of a system of coupled equations. More precisely, subject to observing/measuring growth rate and final epidemic size, we seek to identify parameter values leading to these measurements. We are particularly concerned with disentangling transmission rate from the network density. To do this we define strong and weak identifiability and we find that except for the simplest model, parameters cannot be uniquely determined, that is they are weakly identifiable. This means that there exists multiple solutions (a manifold of infinite measure) which give rise to model output that is close to the data. Identifying, formalising and analytically describing this problem should lead to a better appreciation of the complexity involved in fitting models with many parameters to data.

\end{abstract}


%{\color{red} Re-format 2x2 figures, we're not too interested in results re. $k$ prediction so top-right plot is just unnecessary currently.}

%%%%%%%%%%%%%%%%%%%%%%%%%%%%%%%%%%%%%
\section{Introduction}
%%%%%%%%%%%%%%%%%%%%%%%%%%%%%%%%%%%%%

%\IZK{Choose tau, n pairs that are very different and show that they are very close to the true solution in time.
%Teaser, erdelklodes kelto, es tobb abra nem lesz.}

%\IZK{Model and results in the same section.}

%\IZK{For the compartmental show Determinant nem egyenlo nulla sok beta es gammara.}

%\IZK{DO not plot the intermediary plots and reduce the range around the master plot. Show intersecting curves (inset close to the master point where they intersect) and  temporal solutions.}

Differential-equation-based models are widespread in modelling population dynamics be that in problems arising in ecology, evolution or epidemiology~\cite{anderson1992infectious, blasius2007complex, diekmann2000mathematical,kiss2017mathematics}. Such systems are relatively straightforward to set up and the theory of dynamical systems offers tools to analyse them. Over the past two decades, differential-equation-based models have gained a lot of popularity in modelling epidemics on networks~\cite{Porter2016,kiss2017mathematics}. Such models, often referred to mean-field models, aim to approximate the expected behaviour of some quantities of interest (e.g. expected number of infected individuals in time) and rely on closure assumptions which are needed to produce tractable systems. 
 %Often, network characteristics, such as average degree, degree heterogeneity, appears explicitly in the systems of differential equations, and in expression such as lead eigenvalue, which is essential for stability analysis, and first integrals that allows us to work out quantities such as final epidemic size~\cite{}.
A myriad of ODE-based epidemic models are available~\cite{kiss2017mathematics} with many providing explicit or implicit analytical expressions for quantities such as the basic reproduction number (or leading eigenvalue based on the linear stability analysis around the disease-free steady state), timing and/or peak prevalence, final epidemic size etc. Hence, given a synthetic or real-epidemic and being able to measure a number of the aforementioned quantities, it is of interest to investigate if parameters of the epidemic model that generated the data can be inferred or determined.

\begin{figure}[h!]
     \centering
     \includegraphics[width=0.45\textwidth]{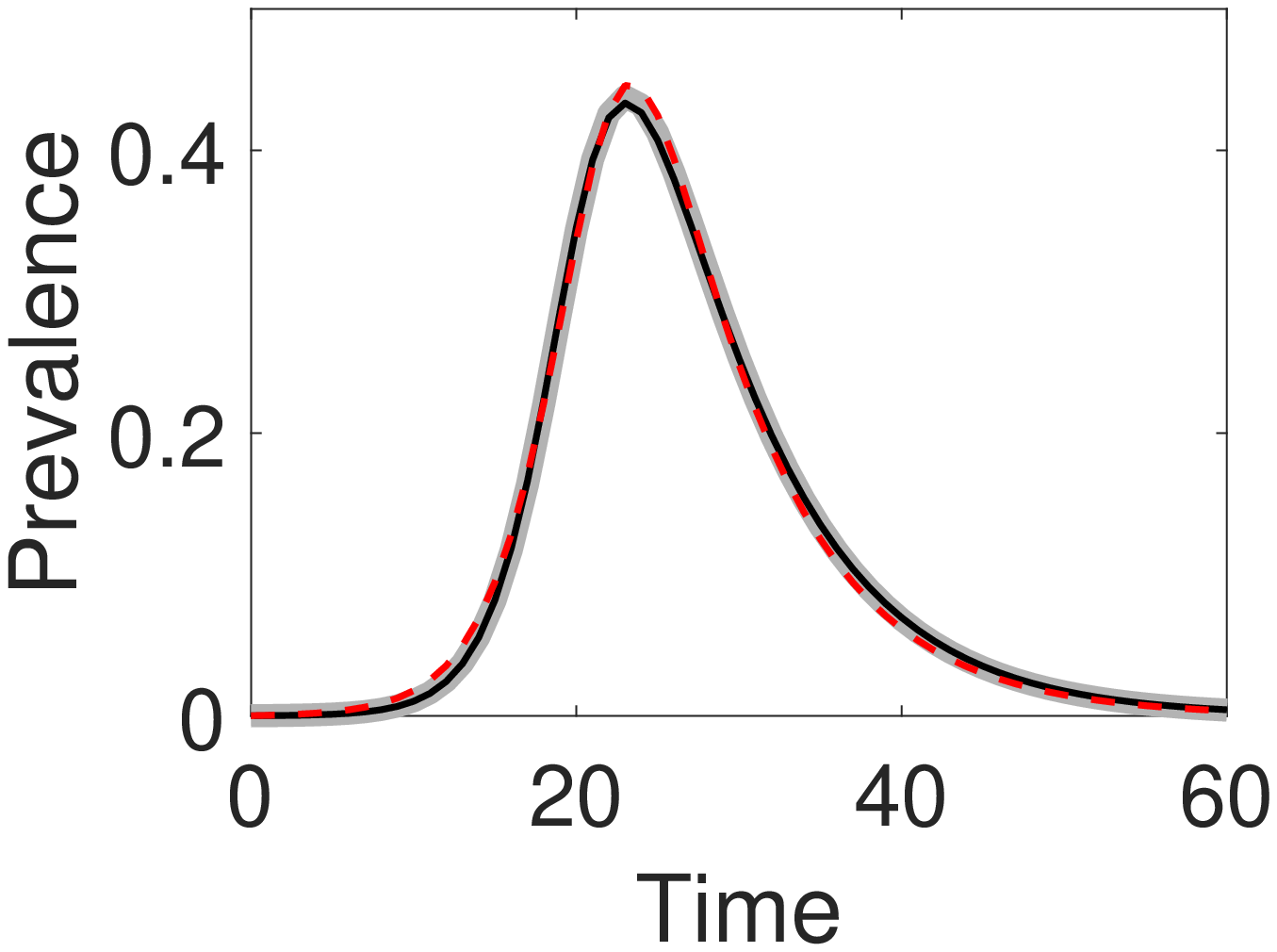}
     \includegraphics[width=0.45\textwidth]{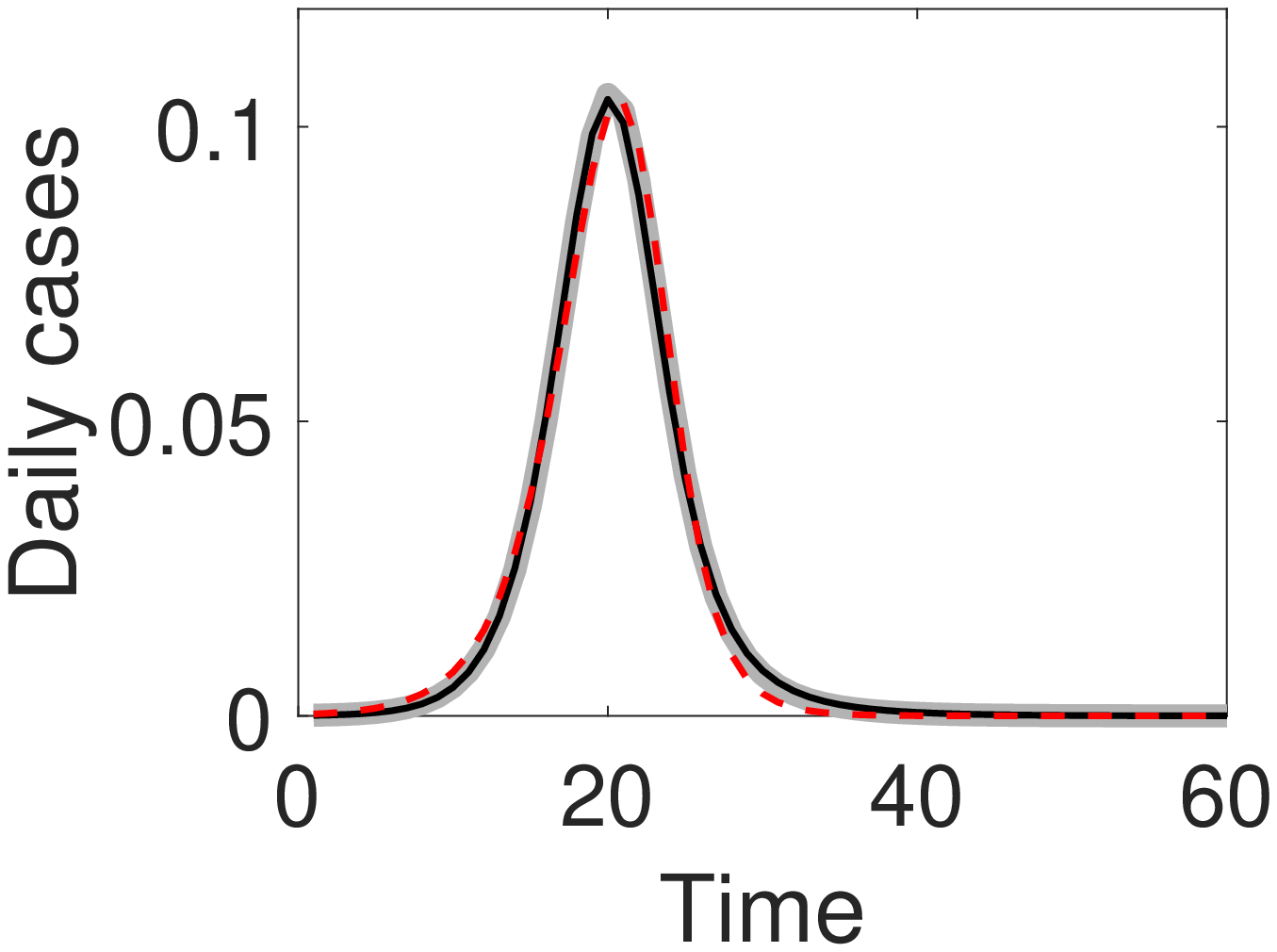}
     \caption{Illustration of how distinct pairs of average degree and transmission rate, $(n, \tau)$ lead to almost indistinguishable time evolution of the prevalence and daily new cases. Baseline values of the parameters are:  average degree $n=6$, $\tau=\gamma R_0/((n-1)-R_0)=0.1429$, with $R_0=2.5$, rate of recovery, $\gamma=1/7$, number of nodes $N=10000$ and epidemic started with one infected individual, with the corresponding output shown by the thick grey lines. The black and red-dashed lines correspond to $(n, \tau)=(8.46, 0.09)$ and $(n, \tau)=(2.454,1.091)$, respectively.}
     \label{fig:two_close_curves}
\end{figure}

Fitting epidemic models to synthetic or real-world data is of great interest as it allows us to infer model parameters which in turn helps us to (i) learn more about the disease, (ii) implement and test control scenarios via simulations, and (iii) make short- or long-term predictions about the epidemic~\cite{chowell2017fitting, King2015}. In many cases, such models can and will be used for parameter estimation and prediction and can suffer of the well-known problem of parameter redundancy and identifiability~\cite{cole2019parameter,villaverde2016structural,gallo2022lack}. This problem has also been highlighted in network-based epidemic models, for example in~\cite{britton2002bayesian}. This problem is not model specific. For example, in Figure~\ref{fig:two_close_curves}, we show that for the pairwise model, Eqs.~\eqref{eq:PWS}-\eqref{eq:PWSS}, it is possible to find distinct sets of parameters whereby the time evolution of prevalence and daily new cases are near indistinguishable. Of course this also implies that the initial growth rate and final epidemic size are also close.

For example, in~\cite{roosa2019assessing}, the authors consider the problem of parameter identifiability in a number of increasingly complex compartmental epidemic models. As the number of states in the model increases so does the number of parameters. While the parameters remained identifiable, in particular the basic reproduction number, the uncertainty around the estimate increased in models with more parameters. On the other hand in~\cite{gallo2022lack}, the authors provide a framework to quantify
how the uncertainty in the data affects the determination of the parameters and the evolution of the unmeasured
variables of a given model. Their approach allows them to characterise different regimes of identifiability and argue that in some cases, such as COVID-19 the lack of identifiability may prevent reliable predictions of the epidemic dynamics. Finally, in~\cite{villaverde2016structural}, the authors argue that structural identifiability in every model should be checked before using the model for inference. But this is seldom done since it involves either complex analytical or numerical calculations. 

%Some models have intrinsic structural identifiability problems which in some cases can be corrected.

In this paper we show that structural unidentifiability is present in a number of network-based epidemic models. This is however, is not due to hidden or unmeasured variables. Moreover, our inference is making use of available analytical formulas for leading eigenvalue, or equivalently growth rate, and final epidemic size. While many of the previous works are concerned with local changes; that is quantifying change in observations induced by a small local change in parameter values, we show that in our models varying the parameters globally leads to small local changes in measurement.   

The paper is structured as follows. In section 2 we describe the general mathematical approach and suggest some definition and ways to formalise the identifiablility problem. In Section 3, we start with simple models such as the well-mixed susceptible-infected-recovered (SIR) compartmental model, followed by more complex models such as, the pairwise, Section 4, and  the edge-based compartmental model. We show that except for the simplest of models, there are clear parameter identifiability problems which we map out and explain analytically, where possible. In models with a larger number of parameters, it is often the case that many different combinations of the model parameters (with many individual parameters being far from their true values) result in output which is consistent with the true epidemic. 
%Section 3 discusses the case of the simplest well-mixed SIR model. The main results are for the pairwise model and edge-based compartmental models in Sections 4 and 5, respectively.
Finally, we provide some discussion and future directions of research.

%From  a parameter inference view point this presents some important challenges which need to be mitigated carefully. We advocate for ...

%%%%%%%%%%%%%%%%%%%%%%%%%%%%%%%%%%%%%
\section{General approach} \label{section2}
%%%%%%%%%%%%%%%%%%%%%%%%%%%%%%%%%%%%%

We are given a system of ODEs involving some parameters:
$$
\dot x(t)= f(x(t),\mu),
$$
where $x(t)\in \mathbb{R}^n$ is the state vector of the system and $\mu\in \mathbb{R}^k$ is the vector of parameters. We observe a derived quantity (e.g., final epidemic size, growth rate) for which data is available. This is given by an observation function $h: \mathbb{R}^n \to \mathbb{R}^m$, i.e. the observation $y$ is
$$
y(t)=h(x(t),\mu) .
$$
\textbf{The goal is to solve the inverse problem}, namely, to determine the parameter $\mu$ based on the observation $y(t)$, note that the observation does not need to be time dependent. This is line with the formulation of a general inverse problem, see ~\cite{cole2019parameter}.

%\url{ https://protect-eu.mimecast.com/s/Y1AxCVQ0YC7rXY8hJC3g4?domain=en.wikipedia.org }

%\url{ https://protect-eu.mimecast.com/s/KS4XCWqjZF27xV8tmtHL8?domain=epubs.siam.org }

\textbf{Our question here is parameter identification}, namely, to understand if it may happen that observations $y$ and $\overline{y}$ corresponding to different parameters $\mu$ and $\overline{\mu}$ are identical or very close to each other.

Thus we are looking for conditions on the parameters ensuring that two observations are close to each other. This can be achieved in two different ways. The first is when the time dependence of the observation is known for all time values, or at least for an observation time-window. The second is when we have formulas for some characteristic quantities of the observation. For example, the derivative of the quantity being observed at the initial instant, $\dot y(0)$, or its limit for large time, $y(\infty)$. These formulas typically involve the unknown parameter values and hence define a system of equations for them. The parameters can be {\it identified} by the model, if this system can be uniquely solved for the parameters. If the opposite holds, then we refer to it as {\it unidentifiability}. That is, we speak about unidentifiability when clearly distinct parameter values satisfy the above system of equations but with some small error.

This can be formulated as follows. Let the system of equations for the parameters take the form $F(\mu)=0$. We call the problem unidentifiable in the strong sense if the system of equations $F(\mu)=0$ has more then one (typically infinitely many) solution. The problem is called \textbf{unidentifiable in the weak sense} if the system $|F(\mu)|<\varepsilon $ is satisfied by a large set of $\mu$ values for any $\varepsilon >0$. In fact, we will show that in our cases the set of $\mu$ values solving $|F(\mu)|<\varepsilon $ has infinite measure. We note that this does not exclude that the equation $F(\mu)=0$ has a unique solution. 

This notion of unidentifiability is related but not identical to the question of sensitive dependence on parameters, which is a notion of local nature. That is expressed in terms of the Jacobian of $F$ at the solution of $F(\mu)=0$. The problem fits into the framework of error analysis and sensitivity analysis that are widely studied important fields of parameter inference. We refer the interested Reader to the books \cite{cacuci2005sensitivity, einarsson2005accuracy}, in which both the introduction to the topic and elaborated examples are available. 

Comparing our definition to that one in \cite{gallo2022lack}, the main novelty in ours is that the inequality $|F(\mu)|<\varepsilon $ holds globally in the parameter space. The definition  in \cite{villaverde2016structural} relates unidentifiability to Kalman's observability condition in an augmented system (where the phase space is extended by new artificial variables representing the parameters). This definition is also local in the parameter space in contrast to our global definition. Another difference between our definition and previous ones is that both \cite{gallo2022lack} and \cite{villaverde2016structural} infers parameters from the time dependence of the solutions while we use exact (not numerical) formulas for some characteristic quantities (leading eigenvalue and final epidemic size). Hence the parameter inference is done by solving a system of equations instead of fitting to time dependent curves. 

We apply this general theory to the SIR compartmental, pairwise and edge-based compartmental models when the observations are the leading eigenvalue at the disease-free steady state  and the final epidemic size. The parameters to be determined, given these observations, are the infection rate $\tau$, the recovery rate $\gamma$ and the average degree of the underlying network, $n$. We will show that for these models (in fact for several other models as well) the leading eigenvalue can be expressed in terms of these parameters as
$$
\lambda = \tau l(n) -\gamma ,
$$
where $l(n)$ is a linear function depending on the model. We will derive an implicit equation for the final epidemic size in each case. It will turn out that this implicit equation contains the parameters $\tau$ and $\gamma$ linearly and $n$ in a non-linear way. The equation can be written in the form
$$
\tau = A(n)\gamma
$$
for all cases, where $A(n)$ is a nonlinear function.

The problem of parameter identification can be formulated as follows. Depending on the choice of the model, i.e. the choice of the functions $l(n)$ and $A(n)$, can the parameters be recovered by solving the two equations above? Since we have two equations for three parameter values, it is obvious that one of the parameters has to be assumed to be given. The recovery rate is more appropriate for being a known value since it depends on epidemiological characteristics. While $\tau$ and especially $n$ are more dependent on the behaviour of the agents and on the network, hence these are more difficult to determine. Our goal will be to solve the above equations for $\tau$ and $n$ with a given value of $\gamma$ and also with given initial conditions. (We note that the initial conditions could also be considered as parameters which makes the problem even more complicated in a real-life situation.)

These equations define two curves in the $(n,\tau)$ parameter plane. The parameter values leading to the desired values of the final size and leading eigenvalue can be obtained as the intersection point of the two curves. The main novelty of the paper is the observation that these curves are very close to each other, hence relatively different parameter values may yield very similar final size and leading eigenvalue. Thus, noisy data may preclude the correct identification of the values of these parameters.

The above system is linear in $\tau$ and $\gamma$ when $n$ is considered to be known. Hence its solvability is easy to check by computing the determinant. For the sake of completeness, this will also be carried out below in each case.

%%%%%%%%%%%%%%%%%%%%%%%%%%%%%%%%%%%%%%%%%%%%%%%%%%%%%
\section{Identifiability in the compartmental SIR model}
%%%%%%%%%%%%%%%%%%%%%%%%%%%%%%%%%%%%%%%%%%%%%%%%%%%%%
The well-known SIR compartmental model takes the form
\begin{align*}
\dot{S} &= -\tau n I\frac{S}{N}, \\
\dot{I} &= +\tau n I\frac{S}{N}-\gamma I .
\end{align*}
Simple differentiation at the disease-free stead state ($S=N$, $I=0$) yields that the leading eigenvalue is $\lambda=\tau n -\gamma$. On the other hand, the final epidemic size is given by the solution of the following implicit equation $R_{\infty}=N-S_0\exp{\left(-\tau n R_{\infty}/N \gamma\right)}$ as it is given in (4.12) in~\cite{kiss2017mathematics}. Let us assume, for sake of simplicity, that $S_0=N$, that is initially there are very few infected and recovered nodes. Then the final size equation can be rearranged to $\tau n r_{\infty} = \gamma \ln (1-r_{\infty})$, where we introduced the fraction $r_{\infty}=R_{\infty}/N$.

Thus the system relating the measured characteristic quantities $\lambda$ and $R_{\infty}$ to the parameters, $\tau$, $\gamma$ and $n$ takes the form
\begin{align}
\tau n-\gamma &=\lambda, \\
\tau n r_{\infty}- \gamma \ln (1-r_{\infty})  &= 0.
\end{align}
This system is linear in $\tau$ and $\gamma$, hence apart from exceptional cases it has a unique solution for $\tau$ and $\gamma$ if $n$ is known and the characteristic quantities of the epidemic, $\lambda$ and $R_{\infty}$, are measured. That is, knowing/measuring the leading eigenvalue and final epidemic size, it is possible to uniquely determine $\tau$ and $\gamma$. However, the parameters $\tau$ and $n$ cannot be obtained from this system, since only their product is determined by the equations. That is, knowing/measuring the leading eigenvalue and final epidemic size, it is not possible to determine the infection rate $\tau$ and average degree $n$. This is the case of strong unidentifiability when the system of equations has infinitely many solutions (if it has a solution at all).

%%%%%%%%%%%%%%%%%%%%%%%%%%%%%%%%%%%%%%%%%%%%%%%
\section{Identifiability in the pairwise SIR model}
%%%%%%%%%%%%%%%%%%%%%%%%%%%%%%%%%%%%%%%%%%%%%%%

The pairwise model focuses on a hierarchical construction where expected number of nodes in state $A$ at time $t$, $[A](t)$,
depends on the expected number of pairs of various types (e.g. $[AB]$) and then, these in turn depend on triples such as $[ABC]$.
Here the counting is done in all possible directions meaning that $[SS]$ pairs are counted twice and and that $[SI]=[IS]$. With this in mind the pairwise model becomes (see e.g. in~\cite{kiss2017mathematics})
\begin{align*}
[\dot{S}] &= -\tau[SI]; \,\,\,[\dot{I}] = \tau[SI]-\gamma[I]; \,\,\, [\dot{R}] = \gamma[I], \\
[\dot{SI}] &= -(\tau+\gamma)[SI] + \tau([SSI]-[ISI]); \,\,\, [\dot{SS}] = -2\tau[SSI].
\end{align*}
This system is not self-consistent as pairs depend on triples and equations for these are needed. To tackle this dependency on higher-order moments the triples in the equation above are closed using the following relation,
\begin{equation*}
[ASB]=\kappa\frac{[AS][SB]}{[S]},
\end{equation*}
where $A, B \in \{A, B\}$. Common choices for $\kappa$ are $(n-1)/n$ and 1. We will consider unidentifyability here for $\kappa=\frac{n-1}{n}$. Applying this closure leads to
\begin{align}
[\dot{S}] &= -\tau[SI], \label{eq:PWS}\\
[\dot{I}] &= \tau[SI]-\gamma[I],  \\
%[\dot{R}] &= \gamma[I], \\
[\dot{SI}] &= -(\tau+\gamma)[SI] + \tau\frac{n-1}{n}\frac{[SI]([SS]-[SI])}{[S]},  \\
[\dot{SS}] &= -2\tau\frac{n-1}{n}\frac{[SS][SI]}{[S]}, \label{eq:PWSS}
\end{align}
which is now a self-contained system.

The leading eigenvalue, resulting from the linear stability analysis around the disease free steady state, ($[S],[I],[SS],[SI]) = (N,0,nN,0)$), can be easily computed from equations~\eqref{eq:PWS}-\eqref{eq:PWSS} as
\begin{equation}
\lambda = \tau(n-2)-\gamma .
\label{eq:PW_LeadEigv}
\end{equation}

An implicit equation for the final number of recovered and susceptible nodes can be derived as it is shown in Section 4.3.4 in~\cite{kiss2017mathematics}. Equation (4.17) there yields the final number of susceptible nodes, $S_{\infty} = N-R_{\infty}$. Let us assume again, for sake of simplicity, that $S_0=N$, that is initially there are very few infected and recovered nodes. Then dividing equation (4.17) in~\cite{kiss2017mathematics} by $N$ and introducing $s_{\infty} = S_{\infty}/N$ leads to
$$
\tau \left( s_{\infty} - s_{\infty}^{2/n} \right)  + \gamma \left( s_{\infty} ^{1/n}-s_{\infty} ^{2/n} \right) = 0.
$$
Thus the system relating the measured characteristic quantities $\lambda$ and $s_{\infty}$ to the parameters, $\tau$, $\gamma$ and $n$ takes the form
\begin{align}
\tau (n-2)-\gamma &=\lambda,  \\
\tau \left( s_{\infty} - s_{\infty}^{2/n} \right)  + \gamma \left( s_{\infty} ^{1/n}-s_{\infty} ^{2/n} \right) &= 0.
\end{align}
This system is linear in $\tau$ and $\gamma$, hence apart from exceptional cases it has a unique solution for $\tau$ and $\gamma$ if $n$ is known and the characteristic quantities of the epidemic, $\lambda$ and $R_{\infty}$ are measured. That is, knowing/measuring the leading eigenvalue and final epidemic size, it is possible to uniquely determine $\tau$ and $\gamma$.

\begin{figure}[h!]
     \centering
     \includegraphics[width=0.95\textwidth]{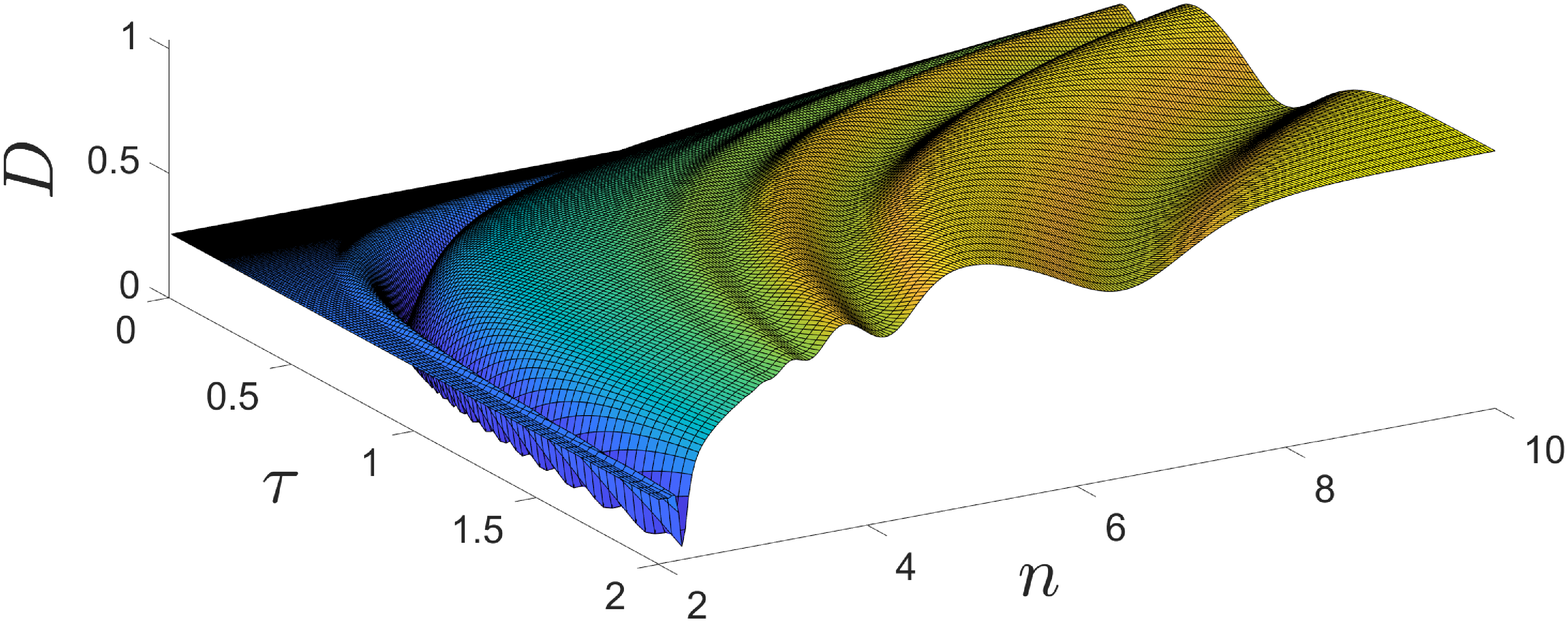}
     \includegraphics[width=0.95\textwidth]{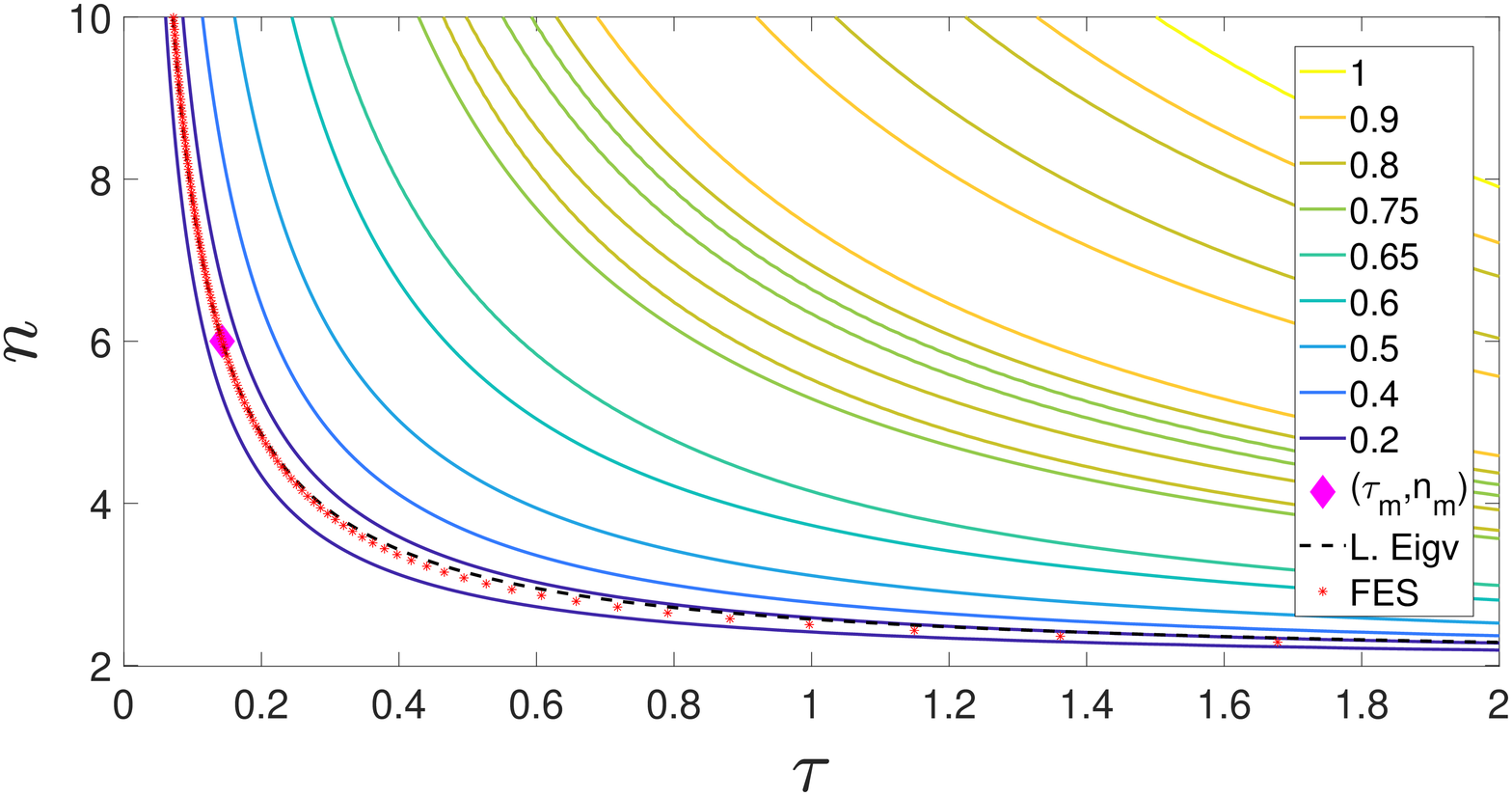}
     \caption{{\bf Top panel:} illustration of the distance profile $D$ between the daily new infections in time for a fixed set of values $(\tau_m,n_m)=(\gamma R_0/((n-1)-R_0)=0.1429,6)$, with $\gamma=1/7$ and $R_0=2.5$, compared to daily new cases for different choices for the values of $(\tau,n)$ pair. Distance measured using an Euclidian norm scaled by the population size $N=10000$. {\bf Bottom panel:} the same as above as contour plot with two additional curves given by the equations for the lead eigenvalue and final epidemic size Eqs.~\eqref{eq:lamPW} and \eqref{eq:sinf}, respectively, where the lead eigenvalue $\lambda$ and $s_{\infty}$ are calculated with $(\tau_m,n_m)$ as given above.}
     \label{fig:3d_and_contour}
\end{figure}

Let us turn now to the parameters $\tau$ and $n$. Now $\gamma$ is considered to be given, and the characteristic quantities of the epidemic, $\lambda$ and $s_{\infty}$ are measured. We can express $\tau$ from the equations above yielding
\begin{align}
    \tau &= \frac{\lambda+\gamma}{n-2}, \label{eq:lamPW} \\
    \tau &= \gamma \frac{s_{\infty} ^{1/n}-s_{\infty} ^{2/n}}{s_{\infty}^{2/n} - s_{\infty}  }  . \label{eq:sinf}
\end{align}
In order to show unidentifiability visually, let us plot the curves given by the above equations in the $(\tau, n)$ plane. We can see in Fig.~\ref{fig:3d_and_contour} (bottom panel) that the two curves are practically indistinguishable. In fact, they have a single intersection point, i.e. the system has a unique solution, but any value of $\tau$ yields a value of $n$ on the hyperbola like curve, that is an approximate solution with high accuracy. 

In fact the experiment that we setup here, and in some of the cases that follow, is that we start with a known set of parameters, often referred to as master set of values. These generate a particular numerical value for the lead eigenvalue, final epidemic size and time evolution of the prevalence or daily new cases. We then ask the questions: are there any other parameter combination ($\tau,n$) that give rise to daily new cases in time that are similar to that obtained by using the master values. The top panel in Figure~\ref{fig:3d_and_contour} shows the euclidean distance between the master daily cases vector and those resulting from ($\tau,n$) pairs chosen between the bounds seen in the figure. 

There are several important features to note about the surface showing the distances. First, there is a clear hyperbola-like valley of minimum points, where any choice of ($\tau,n$) seem to be close enough to the output based on the master values. Several minima are observed which indicate that any kind of optimiser may struggle to find the global optimum. Of course in this thought experiment, there is a unique ($\tau,n$) pair that makes $D=0$. However, given noisy observations, it is easy to see that any values along the hyperbola-like valley may return an acceptable fit, such as the one in Figure~\ref{fig:3d_and_contour}.

The empirical experiment and observations above, can be made more substantial by considering the bottom panel in Figure~\ref{fig:3d_and_contour}. The contour plot is based on the same data as in the surface plot above but with the addition of two curves: that of the leading eigenvalue and final epidemic size, which have unique numerical values determined by the master values and fixed $\gamma$. It is clear that these two curves are indeed close to each other and that they capture the hyperbola-like valley of small values in distance.

\begin{figure}[h!]
     \centering
     \includegraphics[scale=0.5]{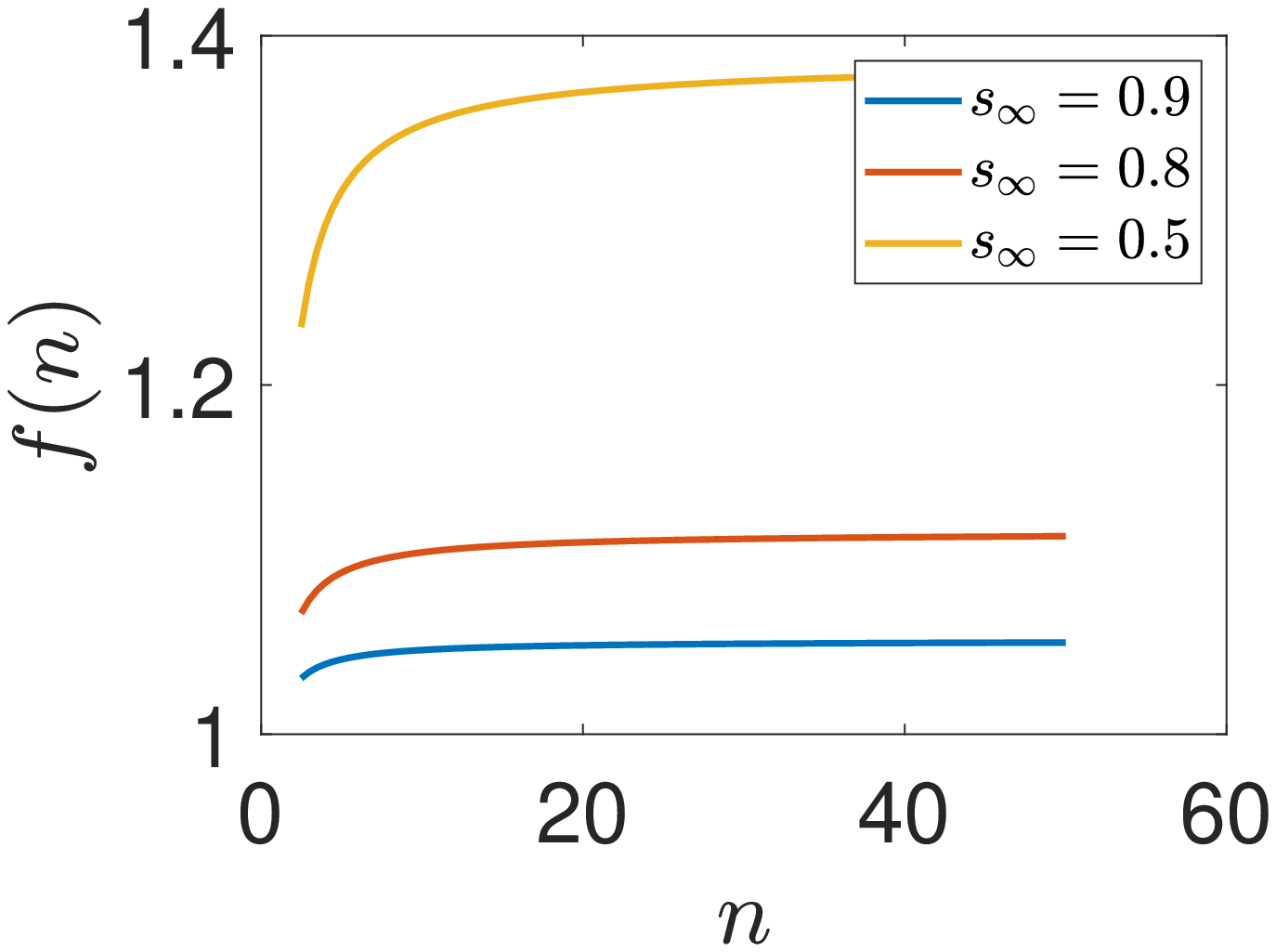}
     \caption{Plots of function $f(n)$. For different values of $s_{\infty}$.}
     \label{fig:f_of_n_sensitivity}
\end{figure}

Beyond this visualisation of unidentifiability, we formally prove it in terms of the definition given in Section \ref{section2}. First, we reduce system \eqref{eq:lamPW}-\eqref{eq:sinf} to a single equation as follows:
$$
\lambda+\gamma = \gamma f(n),
$$
where
$$
f(n)= \frac{s_{\infty} ^{1/n}-s_{\infty} ^{2/n}}{s_{\infty}^{2/n} - s_{\infty}  } (n-2).
$$
We can assume without loss of generality, that the two curves have a common point, i.e. there is a value $n^*$ of $n$ satisfying $\lambda+\gamma = \gamma f(n^*)$. Otherwise, the measurement was so inaccurate that no values of $\tau$ and $n$ could lead to the measured value of $\lambda$ and $s_{\infty}$. Thus the single equation to be solved for the unknown $n$, takes the form
$$
f(n) = f(n^*).
$$
We will prove that this equation does not identify the value of $n$ in the weak sense. In order to do so, we determine the characteristic properties of function $f$. These properties can be easily visualized by plotting the graph of the function for $n>2$, see Figure \ref{fig:f_of_n_sensitivity}. It turns out that the function is very close to a constant, its value changes only slightly from $n=2$ to infinity. For example, in the case $s_{\infty}=0.9$, the functions grows from $1.027$ (at $n=2$) to $1.054$ as $n$ tends to infinity, so the function is constant with accuracy $0.027$. 

Simple application of L'Hospital's rule yields that the limits of $f$ as $n$ tends to $2$ or to infinity exist and their values are
$$
\lim\limits_{n \to 2} f(n) = 2\frac{s_{\infty}-\sqrt{s_{\infty}}}{s_{\infty}\ln{s_{\infty}}} :=f_2
$$
$$
\lim\limits_{n \to \infty} f(n) = \frac{\ln{s_{\infty}}}{s_{\infty}-1} :=f_{\infty}
$$
The next proposition expresses the fact that the measure of the range of this function is small.

\begin{prop}
There exist a number $0<\overline{a}<1$ such that $s_{\infty}>\overline{a}$ implies that $f$ is increasing and $f_2 < f(n) < f_{\infty}$ for all $n>2$. That is the range of $f$ is the interval $(f_2, f_{\infty})$.
\end{prop}

\noindent \textbf{Proof}

Introducing $a=s_{\infty}$, $x=1/n$ and the function
$$
g(x)= \frac{a^{2x} - a^{x} }{a - a^{2x}}\left(\frac{1}{x} -2 \right) ,
$$
we have $f(n)=g(1/n)$, leading to $f'(n)= -g'(1/n)\frac{1}{n^2}$. Hence it is enough to prove that $g'(x)<0$ for all $x\in (0,1/2)$.

Simple differentiation shows that $g'(x)<0$ is equivalent to
$$
\left(\frac{1}{x} -2 \right) \left[ (2a^{2x} - a^{x})(a-a^{2x}) + 2a^{2x}(a^{2x} - a^{x}) \right] \ln a < \frac{1}{x^2} (a^{2x} - a^{x})(a-a^{2x}) ,
$$
that can be rearranged to (by multiplying by $x^2$)
$$
(1-2x) \left[ 2a^{2x+1} -a^{x+1} - a^{3x} \right] \ln a^{x} < (a^{2x} - a^{x})(a-a^{2x}) .
$$
Introducing the new variable $b=a^x$ and returning to $n$ instead of $x$, the desired inequality takes the form (after dividing by $b^3$)
$$
0< n(1-b)(1-b^{n-2})+(n-2)(1+b^{n-2}-2b^{n-1})\ln b := h(b) .
$$
This newly defined function satisfies $h(1)=0$, and elementary differentiation shows that $h'(1)=0=h''(1)$. Moreover, the inequality $h'''(1)<0$ holds. Based on this inequality, it is easy to check that $h$ is positive in a left neighbourhood of $1$, that is there exists a number $\overline{b}<1$, such that $h(b)>0$ holds when $\overline{b}<b<1$. 
%\PLS{Ennek igazolasa, és a h derivaltjainak kiszamolasa még egy harmad oldal lenne, beirjuk ezeket?}\IZK{Ez a te dontesed!}

Let us define the desired number $\overline{a}$ as $\overline{a}=\overline{b}^2$. Then $a>\overline{a}$ is equivalent to $b^n > \overline{b}^2$. On the other hand, $n>2$ and $b<1$ imply that $b^2 >b^n$, hence $b^2 >b^n > \overline{b}^2$, yielding $b> \overline{b}$ leading to $h(b)>0$. This is equivalent to $g'(x)<0$ that we wanted to prove. 

$\Box$

%\PLS{The proof of this proposition could be here, or in an Appendix.}
%\IZK{Maybe in the Appendix?}

We note that numerical evidence shows that the number $\overline{a}$ given by the proposition is relatively small, e.g. $\overline{a}<0.1$. That is, for reasonable values of $s_{\infty}$ the assertions of the proposition hold.

The proposition yields weak unidentifiability as follows. The value of $\gamma$ is considered to be given, and the characteristic quantities of the epidemic, $\lambda$ and $s_{\infty}$ are measured. These determine the unique intersection point $(n^*,\tau^*)$ of the curves given by \eqref{eq:lamPW}-\eqref{eq:sinf}. In other words, $n^*$ is the trivial solution of the reduced single equation $f(n) = f(n^*)$. An approximate solution $n$ satisfies $|f(n) - f(n^*)|<\varepsilon$ with a given positive value of $\varepsilon$. The proposition implies that $|f(n) - f(n^*)|<\varepsilon$ holds for any $n>2$ if $\varepsilon > f_{\infty} - f_2$, which is a small number. An even smaller $\varepsilon$ is achieved if the measured data $\lambda$ and $s_{\infty}$ yield a value of $n^*$ which is larger, i.e. $f(n^*)$ is closer to $f_{\infty}$. Then the value of $\varepsilon$ can be chosen as $\varepsilon = f_{\infty}- f(n^*)$ and then $|f(n) - f(n^*)|<\varepsilon$ holds for $n$ values in a half line, i.e. in a set of measure infinity. This was defined as unidentifiability in the weak sense. 

%\PLS{Here we can show $(\tau,n)$ pairs leading to nearly identical $\lambda$ and $s_{\infty}$, showing again that the measured values of $\lambda$ and $s_{\infty}$ do not determine $\tau$ and $n$. }

%\PLS{On top of that, we can show the full time dependence of the epidemic for these $(\tau,n)$ pairs, verifying that having more information about the epidemic (not only those two characteristic values) does not help to determine $\tau$ and $n$.}

%\IZK{Let's discuss this!}

\begin{figure}[h!]
     \centering
     \includegraphics[scale=0.5]{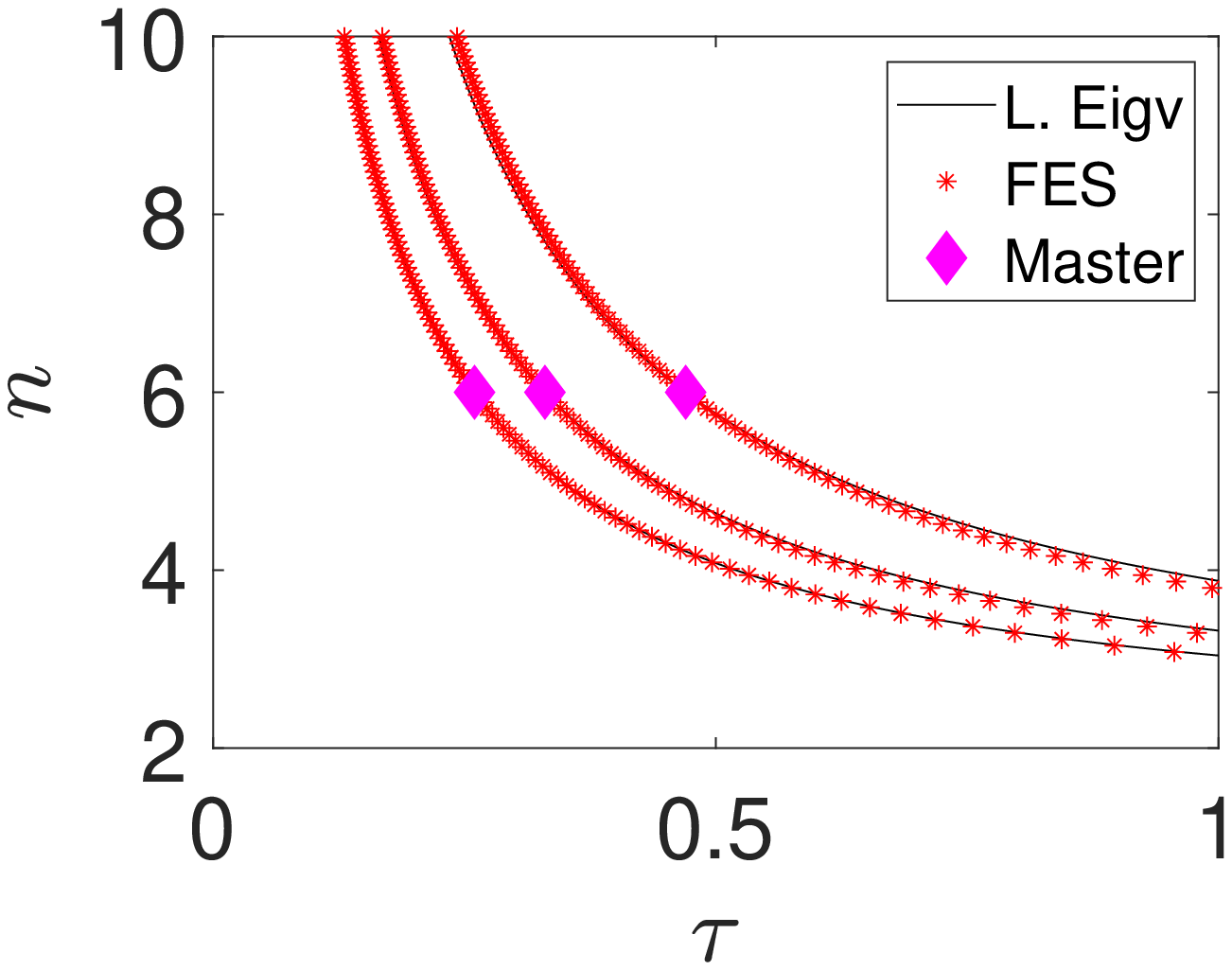}
     \includegraphics[scale=0.5]{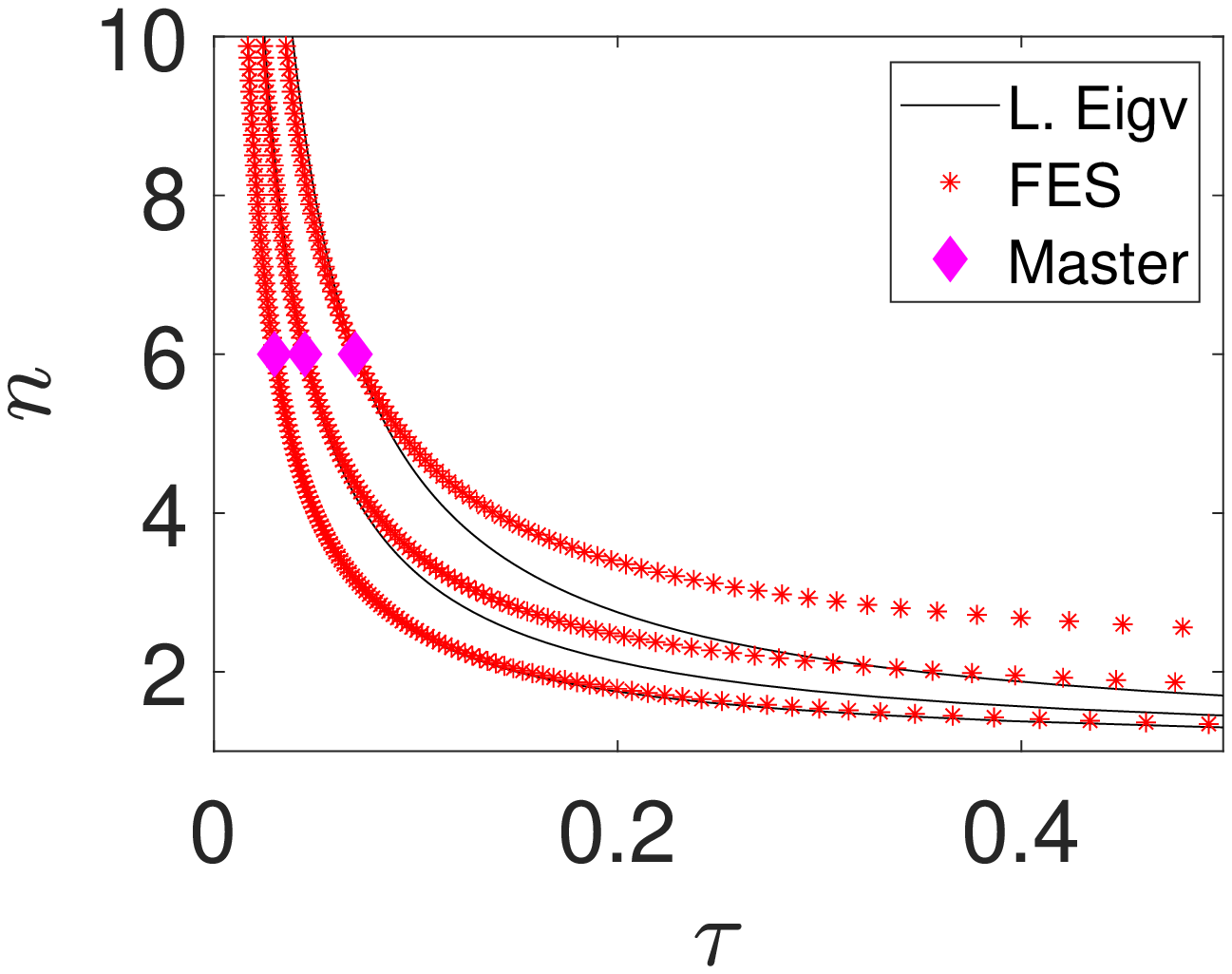}
     \caption{{\bf Left panel:} From left to right curves correspond to solving equations~\eqref{eq:lamPW} and~\eqref{eq:sinf} with the leading eigenvalue and the final epidemic size being set to values obtained by using $\tau=0.26, 0.33, 0.47$, $\gamma=1$ and $n=6$. {\bf Right panel:} Curves given by equations~\eqref{eq:EBCM_UI_Curve_LE} and~\eqref{eq:EBCM_UI_Curve_FS} for values of the transmission rate $\tau=0.03, 0.045, 0.07$ (from left to right). Other parameters are $\gamma=1/7$ and $n=6$. For both plots, the black curve represent ($\tau,n$) pairs where the leading value is that determined by the master values shown as a diamond magenta. Similarly, the red star represent ($\tau, n$) pairs where the final epidemic size is equal to that given by the master values.
     \label{fig:PW_and_EBCM_sensitivity}}
\end{figure}

In Figure~\ref{fig:PW_and_EBCM_sensitivity} we explore the dependency of the weak unidentifiability in the pairwise model on the precise parameters used in the model. The left panel of this figure shows that this feature seem to hold for different parameter combiantions and that we can find infinitely many ($\tau, n$) pairs that lead to a desired eigenvalue and final epidemic size. Moreover, we emphasise again that the two curves do overlap to a great extent and over a large range of parameters.

Before investigating the same problem in a different model, we note that the same calculations for the leading eigenvalue and final epidemic size can be done when the pairwise model is closed with  $\kappa=1$. These calculations lead to

\begin{align}
\tau&=\frac{\lambda_{L}+\gamma}{n-1}, \\ 
 \tau&=\frac{\gamma S_{\infty}(\ln(S_{\infty})-\ln([S](0))}{\frac{[SS](0)}{([S](0))^{2}}\left(S_{\infty}\right)^{2}-S_{\infty}(\ln(S_{\infty})-\ln([S](0))-\left(\frac{[SI](0)}{[S](0)}+\frac{[SS](0)}{[S](0)}\right)S_{\infty}}.
\end{align}
By using the disease-free initial condition, $[S](0)=N$, $[SS](0)=nN$, $[SI](0)=0$ and using that $s_{\infty}=S_{\infty}/N$, the equations above lead to
\begin{align}
\tau&=\frac{\lambda_{L}+\gamma}{n-1}, \\ 
 \tau&=\frac{\gamma \ln(s_{\infty})}{ns_{\infty}-\ln(s_{\infty})-n}.
\end{align}
It urns out that the formulas above are identical to those that we obtain later on for the edge-based compartmental model.

%%%%%%%%%%%%%%%%%%%%%%%%%%%%%%%%%%%%%%%%%%%%%%%
\section{Identifiability in the edge-based compartmental model}
%%%%%%%%%%%%%%%%%%%%%%%%%%%%%%%%%%%%%%%%%%%%%%%
The edge-based compartmental model is given by
\begin{equation}
    \dot{\theta}=-\tau \theta+\tau \phi_{S}(0)\frac{\psi'(\theta)}{\psi'(1)}+\gamma(1-\theta)+\tau \phi_{R}(0)=f(\theta),
    \label{eq:EBCM}
\end{equation}
where $\theta$ denotes the probability that a random neighbour $\nu$ of a random, initially susceptible test node $u$ has not yet passed infection to $u$. Furthermore, $\phi_{S}(0)$ and $\phi_{R}(0)$ are the probabilities that, at $t=0$, the random neighbour $\nu$ of a random, initially susceptible test node $u$ is susceptible and recovered, respectively. Typical initial condition for this system are: $\theta(0)=1$, $\phi_{R}(0)=0$, $\theta(0)=1$ and $\phi_{S}(0)=1-\varepsilon$.

We now consider the case of $\psi(x)=\exp(n(x-1))$; that is a network with Poisson degree distribution with mean $n$. Linearising around $\theta=1$, we obtain

\begin{equation}
    f'(\theta){\Large{|}}_{\theta=1}=\tau \phi_{S}(0)\frac{\psi''(1)}{\psi'(1)}-\tau-\gamma=\tau \frac{n^2}{n}-\tau-\gamma=\tau(n-1)-\gamma=\lambda_{L}^{EBCM},
    \label{eq:lam_EBCM}
\end{equation}

The final epidemic size can also be worked out by finding $\lim_{t\rightarrow \infty}\theta(t)=\theta_{\infty}$ and using that the final proportion of susceptible left in the population is $s_{\infty}=\psi(\theta_{\infty})$. Setting the right hand side of equation~\eqref{eq:EBCM} to zero, an implicit equation for $\theta_{\infty}$ follows,

\begin{equation}
    (\tau+\gamma)\theta_{\infty}-\gamma-\tau \psi_{S}(0)e^{n(\theta_{\infty}-1)}=0.
    \label{eq:FES_EBCM_THETA}
\end{equation}

Since $s_{\infty}=\exp(n(\theta_{\infty}-1))$, equation~\eqref{eq:FES_EBCM_THETA} can be recast in terms of $s_{\infty}$ and yields
\begin{equation}
   \tau n\psi_{S}(0)s_{\infty} -(\tau+\gamma)\ln(s_{\infty})-\tau n=0.
    \label{eq:FES_EBCM_SINFTY}
\end{equation}
We are now in a position to write down a system of equations based on ~\eqref{eq:lam_EBCM} and ~\eqref{eq:FES_EBCM_SINFTY}

\begin{align}
\tau&=\frac{\lambda +\gamma}{n-1}, \label{eq:EBCM_UI_Curve_LE}\\
 \tau&=\frac{\gamma \ln(s_{\infty})}{n\psi_{S}(0)s_{\infty}-\ln(s_{\infty})-n}
% (gam_master*fs_master*(log(fs_master)-log(S0)))/...
%        (  (SS0/(S0^2))*fs_master^2 ...
%        -fs_master*(log(fs_master)-log(S0))...
%        -(SI0/S0)*(fs_master)...
%        -(SS0/S0)*(fs_master));
\label{eq:EBCM_UI_Curve_FS}
\end{align}

%\begin{figure}[h!]
%     \centering
%     \includegraphics[scale=0.75]{Figures/Fig3_EBCM_param_sensitivity.eps}
%     \caption{Illustration of equations \eqref{eq:EBCM_UI_Curve_System} with $(\tau_m,n_m)=(6,0.04)$ (diamond marker), $\gamma=1/7$. In this case $s_{\infty}^{*}\simeq0.56$ and $\lambda_{L}^{EBCM,*}\simeq0.057$}
  %   \label{fig:3d_and_contour}
%\end{figure}

%\textbf{Itt mar nem csinalnek formalis bizonyitast epsilonnal. Csak ket abra eleg lenne. Az egyik, hogy ez a ket hiperbola jellegu gorbe kozel van egymashoz. A masik az idofugges ket (tau,n) parnal, hogy mennyire egyforma.}

Thes curves are shown in the right panel of Figure \ref{fig:PW_and_EBCM_sensitivity}. It can be seen that the two curves are close to each other. The coincidence is more emphasised when $s_{\infty}$ is larger, i.e. the final epidemic size is smaller. 

Beyond this visualization of unidentifiability, we formally prove that in terms of the definition given in Section \ref{section2}. First, we reduce the above system to a single equation as follows:
$$
\lambda+\gamma = \gamma q f(n),
$$
where
$$
f(n)= \frac{n-1}{n-q} , \quad \mbox{ and } \quad q=\frac{\ln s_{\infty}}{s_{\infty} - 1} >1.
$$
We can assume without loss of generality, that the two curves have a common point, i.e. there is a value $n^*$ of $n$ satisfying $\lambda+\gamma = \gamma q f(n^*)$. Otherwise, the measurement was so inaccurate that no values of $\tau$ and $n$ could lead to the measured value of $\lambda$ and $s_{\infty}$. Thus the single equation to be solved for the unknown $n$, takes the form
$$
f(n) = f(n^*).
$$
This equation does not identify the value of $n$ in the weak sense. By plotting the graph of $f$, it turns out that the function is very close to a constant, its value changes only slightly from large values of $n$ to infinity. For example, in the case $s_{\infty}=0.9$, the functions changes from $1.006$ (at $n=10$) to $1$ as $n$ tends to infinity, so the function is constant with accuracy $0.006$ in the infinite half-line $n>10$. In general, one can directly see that $f$ is decreasing and its limit is $1$ as $n$ tends to infinity. Similarly to the case of the pairwise model, weak unidentifiability follows from the fact that the function f is to a constant.

%\newpage

%%%%%%%%%%%%%%%%%%%%%%%%%%%%%%%%%%%%%%%%%%%%%%%%%%%%
\section{Discussion}
%%%%%%%%%%%%%%%%%%%%%%%%%%%%%%%%%%%%%%%%%%%%%%%%%%%%
In this paper we study the identifiability of parameters in network-based epidemic models. We find that network density and the transmission rate cannot be disentangled. More formally this means that when considering these parameters, the model is structurally not identifiable. Preliminary analysis suggests that combinations of $n$ and $\tau$ and other parameters are better behaved, for example when packaged into the expression for $R_0$; this is in line with how to deal with identifiability problems~\cite{villaverde2016structural}. 

In an ideal situation the leading eigenvalue and final epidemic size can be measured to any desired accuracy. Assuming that this is the case, an exhaustive search in the parameter space, again to arbitrary precision, would be able to identify the precise parameters which generated the data. However, real-life observations are noisy and even a small measurement error can lead to a significant shift in the values of the inferred parameters. This leads to what we call weak unidentifiability. 

Contact patterns and the transmission of the disease across a link are strongly related and often are difficult to disentangle. Intuitively, it is known that dense networks with low transmission rate and spare networks with high transmission rate can produce similar epidemics. In fact, our hyperbolas trace out and connect these regimes.
Of course, in this case a Bayesian approach may alleviate the problem. With more and more mobility data becoming available as well as data from contact surveys means that contact networks can be characterised sufficiently in order to produce meaningful estimates from complex models.

In terms of future work, we believe that there is value in carrying out a systematic search over the parameter space to identify areas where the unidentifiability is the most significant. Our preliminary analysis shows that this is both model and parameter dependent. We also note that unidentifiability seems to be more marked for less severe epidemics. For larger epidemics, the overlap between the two hyperbolas decreases, meaning that parameters are easier to identify. 

\newpage

%%%%%%%%%%%%%%%%%%%%%%%%%%%%%%%%%%%%%%%%%%%%%%%%%%%%
\section*{Acknowledgements}
%%%%%%%%%%%%%%%%%%%%%%%%%%%%%%%%%%%%%%%%%%%%%%%%%%%%
Istv\'an Z. Kiss  acknowledge support from the Leverhulme Trust for the Research Project Grant RPG-2017-370. P\'eter L. Simon acknowledges support from the Hungarian Scientific Research Fund, OTKA, (grant no. 135241) and from the Ministry of Innovation and Technology NRDI Office within the framework of the Artificial Intelligence National Laboratory Program.

\section*{Data availability} The datasets generated during and/or analysed during the current study are available from the corresponding author on reasonable request.

\section*{Competing interests} The authors declare no competing interests.

\newpage

\bibliography{project1}

\begin{thebibliography}{10}

\bibitem{anderson1992infectious}
Roy~M Anderson and Robert~M May.
\newblock {\em Infectious diseases of humans: dynamics and control}.
\newblock Oxford university press, 1992.

\bibitem{blasius2007complex}
Bernd Blasius, J~rgen Kurths, and Lewi Stone.
\newblock {\em Complex Population Dynamics: Nonlinear Modeling in Ecology,
  Epidemiology, and Genetics}, volume~7.
\newblock World Scientific, 2007.

\bibitem{britton2002bayesian}
Tom Britton and Philip~D O'Neill.
\newblock Bayesian inference for stochastic epidemics in populations with
  random social structure.
\newblock {\em Scandinavian Journal of Statistics}, 29(3):375--390, 2002.

\bibitem{cacuci2005sensitivity}
Dan~G Cacuci, Mihaela Ionescu-Bujor, and Ionel~Michael Navon.
\newblock {\em Sensitivity and uncertainty analysis, volume II: applications to
  large-scale systems}.
\newblock CRC press, 2005.

\bibitem{chowell2017fitting}
Gerardo Chowell.
\newblock Fitting dynamic models to epidemic outbreaks with quantified
  uncertainty: A primer for parameter uncertainty, identifiability, and
  forecasts.
\newblock {\em Infectious Disease Modelling}, 2(3):379--398, 2017.

\bibitem{cole2019parameter}
Diana~J Cole.
\newblock Parameter redundancy and identifiability in hidden markov models.
\newblock {\em Metron}, 77(2):105--118, 2019.

\bibitem{diekmann2000mathematical}
Odo Diekmann and Johan Andre~Peter Heesterbeek.
\newblock {\em Mathematical epidemiology of infectious diseases: model
  building, analysis and interpretation}, volume~5.
\newblock John Wiley \& Sons, 2000.

\bibitem{einarsson2005accuracy}
Bo~Einarsson.
\newblock {\em Accuracy and reliability in scientific computing}.
\newblock SIAM, 2005.

\bibitem{gallo2022lack}
Luca Gallo, Mattia Frasca, Vito Latora, and Giovanni Russo.
\newblock Lack of practical identifiability may hamper reliable predictions in
  covid-19 epidemic models.
\newblock {\em Science advances}, 8(3):eabg5234, 2022.

\bibitem{King2015}
Aaron~A. King, Matthieu~Domenech {De Cell{\'{e}}s}, Felicia~M.G. Magpantay, and
  Pejman Rohani.
\newblock {Avoidable errors in the modelling of outbreaks of emerging
  pathogens, with special reference to Ebola}.
\newblock {\em Proceedings of the Royal Society B: Biological Sciences},
  282(1806):0--6, 2015.

\bibitem{kiss2017mathematics}
Istv{\'a}n~Z Kiss, Joel~C Miller, P{\'e}ter~L Simon, et~al.
\newblock Mathematics of epidemics on networks.
\newblock {\em Cham: Springer}, 598:31, 2017.

\bibitem{Porter2016}
Mason~A. Porter and James~P. Gleeson.
\newblock {Dynamical Systems on Networks}.
\newblock pages 49--51. 2016.

\bibitem{roosa2019assessing}
Kimberlyn Roosa and Gerardo Chowell.
\newblock Assessing parameter identifiability in compartmental dynamic models
  using a computational approach: application to infectious disease
  transmission models.
\newblock {\em Theoretical Biology and Medical Modelling}, 16(1):1--15, 2019.

\bibitem{villaverde2016structural}
Alejandro~F Villaverde, Antonio Barreiro, and Antonis Papachristodoulou.
\newblock Structural identifiability of dynamic systems biology models.
\newblock {\em PLoS computational biology}, 12(10):e1005153, 2016.

\end{thebibliography}
\bibliographystyle{plain}
\newpage

%%%%%%%%%%%%%%%%%%%%%%%%%%%%%%%%%%%%%%%%%%%%%%%%%%%%
%\section{Appendix}
%%%%%%%%%%%%%%%%%%%%%%%%%%%%%%%%%%%%%%%%%%%%%%%%%%%%

\end{document}